# A limited in bandwidth uniformity for the functional limit law of the increments of the empirical process

**Davit Varron**

*Département de Mathématiques, Université de Fracche-Comté, 16 Route de Gray, 25000 Besançon, France*
*e-mail:* dvarron@univ-fcomte.fr

**Abstract:** Consider the following local empirical process indexed by $K \in \mathcal{G}$, for fixed $h > 0$ and $z \in \mathbb{R}^d$:

$$G_n(K,h,z) := \sum_{i=1}^{n} K\Big(\frac{Z_i - z}{h^{1/d}}\Big) - \mathbb{E}\Big(K\Big(\frac{Z_i - z}{h^{1/d}}\Big)\Big),$$

where the $Z_i$ are i.i.d. on $\mathbb{R}^d$. We provide an extension of a result of Mason (2004). Namely, under mild conditions on $\mathcal{G}$ and on the law of $Z_1$, we establish a uniform functional limit law for the collections of processes $\Big\{G_n(\cdot, h_n, z),\ z \in H,\ h \in [h_n, \mathfrak{h}_n]\Big\}$, where $H \subset \mathbb{R}^d$ is a compact set with nonempty interior and where $h_n$ and $\mathfrak{h}_n$ satisfy the Csörgő-Révész-Stute conditions.

**AMS 2000 subject classifications:** Primary 60F15, 60F17; secondary 62G07.
**Keywords and phrases:** Empricial processes, Functional limit theorems, Strong theorems, Density estimation, $\LaTeX\ 2_\varepsilon$.

## Contents

## 1. Introduction and statement of main results

### 1.1. Introduction

Let $(Z_i)_{i\geq 1}$ be a independent, identically distributed sample taking values in $\mathbb{R}^d$. Since the pioneering works of Stute (1982), several researchers have investigated the limit behaviour of the functional increments of the empirical process, which are defined as follows, for fixed $h > 0$ and $z \in \mathbb{R}^d$:

$$\Delta\alpha_n(\cdot, h, z) : s \mapsto n^{-1/2}\left(\sum_{i=1}^n 1_{[z, z+h_n^{1/d}s]}(Z_i) - \mathbb{P}\Big(Z_1 \in [z, z+h_n^{1/d}s]\Big)\right), s \in [0,1]^d.$$

Here we write $[a, b] := [a_1, b_1] \times \cdots \times [a_d, b_d]$ for $a, b \in \mathbb{R}^d$. Deheuvels and Mason (1992) have provided a uniform functional limit law (UFLL) for the following collections of functional increments:

$$\widetilde{\Theta}_n := \left\{\frac{1}{(2h_n \log(1/h_n))^{1/2}} \Delta\alpha_n(\cdot, h_n, z),\ z \in [0, 1 - h_n]\right\}, \tag{1}$$

when $d = 1$, $Z_1$ is uniform on $[0, 1]$ and $h_n$ satisfies the Csörgő-Révész-Stute (CRS) conditions (see (HV1)–(HV3) below). Implicit in their result is the UFLL for $\widetilde{\Theta}_n$ when $(Z_i)_{i\geq 1}$ take values in $\mathbb{R}$ and have a density $f$ that is continuous and strictly positive on an open set $O$, and when $z$ appearing in (1) ranges in a bounded interval $H \subset O$. However, the extension of this result to the multivariate case (d¿1) remained an open problem for almost a decade. Recently, Mason (2004) (see also Einmahl and Mason (2000)) solved this problem by combining the techniques of Deheuvels and Mason (1992) with recent tools in general empirical process theory. Namely, he obtained asymptotic results in a more general framework, considering the following type of stochastic processes indexed by $K$:

$$G_n(K, h, z) = \sum_{i=1}^n \left[K\left(\frac{Z_i - z}{h^{1/d}}\right) - \mathbb{E}\left(K\left(\frac{Z_i - z}{h^{1/d}}\right)\right)\right]. \tag{2}$$

Here, $K$ ranges through a class of functions $\mathcal{G}$ satisfying some conditions that will stated later (see (HK1)–(HK5) in the sequel). More precisely, Mason established a UFLL for the following the sets of processes, as $n \to \infty$,

$$\widetilde{\Theta'_n} := \left\{\frac{G_n(\cdot, h_n, z)}{\sqrt{2f(z)nh_n \log(1/h_n)}},\ z \in H\right\},$$

where $H$ is a compact set of $\mathbb{R}^d$ with nonempty interior. To cite his result, we have to recall the basic assumptions made in Mason (2004).

We say that a sequence of constants satisfies the Csörgő-Révész-Stute (CRS) conditions whenever

(HV1) $h_n \downarrow 0,\ 0 < h_n < 1,\ nh_n \uparrow \infty,$
(HV2) $\lim_{n\to\infty} nh_n/\log n = \infty,$
(HV3) $\lim_{n\to\infty} \log(1/h_n)/\log\log n = \infty.$

Let $\mathcal{G}$ be a class of real Borel functions on $\mathbb{R}^d$. Set $I^d := [0,1]^d$ and

$$\mathcal{F} := \Big\{ K(\lambda(\cdot - z)),\ z \in \mathbb{R}^d,\ \lambda > 0,\ K \in \mathcal{G} \Big\}.$$

Let $||\cdot||_{\mathbb{R}^d}$ be the euclidian norm on $\mathbb{R}^d$. We make the following assumptions on $\mathcal{G}$.

(HK1) $\lim_{||u||_{\mathbb{R}^d}\to 0} \sup_{K\in\mathcal{G}} \int_{\mathbb{R}^d} (K(x) - K(x+u))^2 dx = 0,$
$\lim_{\lambda\to 1} \sup_{K\in\mathcal{G}} \int_{\mathbb{R}^d} \big(K(\lambda x) - K(x)\big)^2 dx = 0,$
(HK2) $\forall K \in \mathcal{G},\ \sup_{x\in\mathbb{R}^d} \mid K(x) \mid \leq 1,$
(HK3) $\forall K \in \mathcal{G},\ \forall x \notin I^d,\ K(x) = 0,$
(HK4) $\exists C_0 > 0,\ v_0 > 0,\ \forall \epsilon \in (0,1),\ \mathcal{N}(\epsilon, \mathcal{F}) \leq C_0 \epsilon^{-v_0}.$

Here, $\mathcal{N}(\epsilon, \mathcal{F})$ denotes the uniform covering number of $\mathcal{F}$ for $\epsilon$ and the class of norms $\{L_2(\mathbb{P})\}$, with $\mathbb{P}$ varying in the set of all probability measures on $\mathbb{R}^d$, and taking the $F \equiv 1$ as an envelope function for the class $\mathcal{F}$ (for more details, see Van der Vaart and Wellner (1996), p. 83–84, with $r = 2$). To overcome any measurability problem, we make the following assumption.

$(HK5)$ $\mathcal{F}$ is pointwise separable (see Van der Vaart and Wellner (1996), p. 110–111).

Let $L_2^*(\mathcal{G})$ be the Hilbert subspace of $L_2(\mathbb{R}^d, \lambda)$ spanned by $\mathcal{G}$, where $\lambda$ denotes the Lebesgue measure on $\mathbb{R}^d$. The (rate) function $J$ that rules the large deviation properties of the isonormal Gaussian process generated by $\big(L_2^*(\mathcal{G}), \lambda\big)$ can be defined, for a function $\Psi : \mathcal{G} \mapsto \mathbb{R}$, by

$$J(\Psi) := \inf\left\{ \int_{\mathbb{R}^d} g^2 d\lambda,\ g \in L_2^*(\mathcal{G}),\ \forall K \in \mathcal{G},\ \Psi(K) = \int_{\mathbb{R}^d} gK d\lambda \right\},$$

with the implicit convention $\inf\ \emptyset = \infty$. Now set $\mathbb{K} = \mathbb{K}_{\mathcal{G}} := \{\Psi :\ J(\Psi) \leq 1\}$. Let $\mathcal{B}(\mathcal{G})$ be the set of all real bounded functions on $\mathcal{G}$. Denote by $|\cdot|_d$ the usual max norm on $\mathbb{R}^d$, namely

$$\mid z \mid_d := \max_{j=1,\ldots,d} \mid z_j \mid . \tag{3}$$

We make a last assumption upon the law of the $Z_i$ (recall that $H$ is a compact set with nonempty interior).

(H f) There exists $\alpha$ such that $Z_1$ has a density $f$ that is continuous and strictly positive on the set

$$H^\alpha := \Big\{ z \in \mathbb{R}^d : \inf_{y \in H} \mid z - y \mid_d < \alpha \Big\}. \tag{4}$$

Under all the above assumptions, Mason established the following result.

**Theorem** (Mason, 2004)**.** *Let $H$ be a compact subset of $\mathbb{R}^d$ with nonempty interior. Let $(Z_i)_{i\geq 1}$ be an i.i.d. sequence of random variables satisfying (H f). Let $(h_n)_{n\geq 1}$ be a sequence of constants fulfilling (HV1)–(HV3) and let $\mathcal{G}$ be a class of real Borel functions satisfying (HK1)–(HK5). Then we have almost surely:*

$$(i)\ \lim_{n\to\infty} \sup_{z\in H} \inf_{\Psi\in\mathbb{K}} \left|\left| \frac{G_n(\cdot, z, h_n)}{\big(2f(z)nh_n \log(1/h_n)\big)^{1/2}} - \Psi \right|\right|_{\mathcal{G}} = 0,$$

$$(ii)\ \forall \Psi \in \mathbb{K},\ \lim_{n\to\infty} \inf_{z\in H} \left|\left| \frac{G_n(\cdot, z, h_n)}{\big(2f(z)nh_n \log(1/h_n)\big)^{1/2}} - \Psi \right|\right|_{\mathcal{G}} = 0.$$

The author proved this result by combining the ideas of Einmahl and Mason (2000) with some recent results in large deviation theory (see Arcones (2003, 200)), Gaussian approximation results for finite dimensional laws (Zaitsev (1987a,b)) and sharp bounds for empirical processes (see Talagrand (1994), or Bousquet (2002) and Klein (2002) for sharper bounds). In the present paper, we show that the arguments of Mason can be efficiently used to enrich his theorem with an additional uniformity in $h \in [h_n, \mathfrak{h}_n]$, under some mild conditions. The remainder of this article is organised as follows. The main result is given by Theorem 1 in §1.2. The proof follows in §2 and is divided into two parts. In §2.1 we sate a large deviation result and a concentration inequality. These two results are somewhat straightforward in regard to the works of Arcones (2003) and Einmahl and Mason (2000). They will play a crucial role in our proof of Theorem 1, which is written in §2.2.

### *1.2. Statement of the result*

We provide in the present paper an extension of the just mentioned theorem of Mason (2004) showing that his UFLL still holds uniformly in $h_n \leq h \leq \mathfrak{h}_n$, provided that both $(h_n)_{n\geq 1}$ and $(\mathfrak{h}_n)_{n\geq 1}$ satisfy (HV1)–(HV3).

**Theorem 1.** *Let $H$ be a compact subset of $\mathbb{R}^d$ with nonempty interior. Let $(Z_i)_{i\geq 1}$ be an i.i.d. sequence of random variables satisfying (H f). Let $(h_n)_{n\geq 1}$ and $(\mathfrak{h}_n)_{n\geq 1}$ be two sequences of positive numbers satisfying (HV1)–(HV3) as well as $\mathfrak{h}_n > 2h_n$. Then we have almost surely:*

$$(i)\ \lim_{n\to\infty} \sup_{h_n\leq h\leq \mathfrak{h}_n, z\in H} \inf_{\Psi\in\mathbb{K}} \left|\left| \frac{G_n(\cdot, h, z)}{\sqrt{2f(z)nh\log(1/h)}} - \Psi \right|\right|_{\mathcal{G}} = 0,$$

$$(ii)\ \forall \Psi \in \mathbb{K},\ \lim_{n\to\infty} \sup_{h_n\leq h\leq \mathfrak{h}_n} \inf_{z\in H} \left|\left| \frac{G_n(\cdot, h, z)}{\sqrt{2f(z)nh\log(1/h)}} - \Psi \right|\right|_{\mathcal{G}} = 0.$$

**Sketch our proof**: Roughly speaking, our proof is divided into the following steps:

- The proof of Mason can be very crudely summed up as follows: given properly chosen sequences of events $(E_n(\epsilon, h_n))_{n\geq 1}$, he proves that, for fixed $\epsilon > 0$ we have, for all large $n$

$$\mathbb{P}\big(E_n(\epsilon, h_n)\big) \leq h_n^{\delta}, \tag{5}$$

for some $\delta > 0$. Then he makes use of the fact that $(h_n)_{n\geq 1}$ satisfies conditions (HV1)–(HV3) to achieve his goals, by making use of usual blocking techniques along with the Borel-Cantelli lemma.
- Given $\rho > 1$, we discretise $[h_n, \mathfrak{h}_n]$ into the following grid of size $R_n \approx \log(\mathfrak{h}_n/h_n)/\log(\rho)$:

$$\{h_{n,0}, h_{n,1}, h_{n,2}\ldots, h_{n,R_n}\} = \{h_n, \rho h_n, \rho^2 h_n, \ldots, \mathfrak{h}_n\}. \tag{6}$$

- For fixed $\epsilon$, we show that $\mathbb{P}(E_n(\epsilon, h_{n,l})) \leq h_{n,l}^{\delta}$ *uniformly* in $l$, for some $\delta > 0$. To do this, we make use of argument that are very similar to those of Mason for proving (5), but taking additional care to get inequalities uniformly in $l$ (see our *key argument* in §2.2, Step 1). Indeed, we had to write an concentration inequality (see Proposition 2.2), which is somehow a finite distance version of the inequality used by Mason.
- Then we write, $\Delta_n$ denoting a term of oscillation of our proceses between two consecutives $h_{n,l}$,

$$\mathbb{P}\Big(\bigcup_{h\in[h_n,\mathfrak{h}_n]} E_n(\epsilon, h)\Big) \leq \mathbb{P}\Big(\bigcup_{l=0}^{R_n} E_n(\epsilon, h_{n,l})\Big) + \Delta_n$$
$$\leq \sum_{l=0}^{R_n} h_{n,l}^{\delta} + \Delta_n \leq h_n^{\epsilon} \sum_{l=0}^{R_n} \rho^{l\delta} + \Delta_n \simeq \frac{\mathfrak{h}_n^{\delta}}{\rho^{\delta} - 1} + \Delta_n.$$

Thus, we can make use of the fact that $(\mathfrak{h}_n)_{n\geq 1}$ satisfies (HV1)–(HV3) and continue our proof as in the proof of Mason. The oscillation term $\Delta_n$ is controlled by the concentration inequality of Proposition 2.2. We now focus on some corollaries of Theorem 1. Denote by $g_{n,h,z}$ a non usual form of functional increments of the empirical process, namely:

$$g_{n,h,z}(s) := \frac{1}{nh}\sum_{i=1}^{n} 1_{[s,1]}\Big(\frac{Z_i - z}{h^{1/d}}\Big) - \mathbb{E}\Big(1_{[s,1]}\Big(\frac{Z_i - z}{h^{1/d}}\Big)\Big),\ s \in [0,1]^d, \tag{7}$$

with the notation $[s,1] := [s_1,1]\times\cdots\times[s_d,1],\ s = (s_1,\ldots,s_d) \in [0,1]^d$. Applying Theorem 1 respectively to the particular class of indicator functions of the sets $[s,1],\ s \in [0,1]^d$, we obtain, almost surely:

$$\lim_{n\to\infty} \sup_{h\in[h_n,\mathfrak{h}_n]} \sup_{z\in H} \inf_{g\in\mathcal{S}} \sqrt{\frac{nh}{2f(z)\log(1/h)}}\ ||\ g_{n,h,z} - g\ ||_{[0,1]^d} = 0, \tag{8}$$

$$\forall g \in \mathcal{S},\ \lim_{n\to\infty} \sup_{h\in[h_n,\mathfrak{h}_n]} \inf_{z\in H} \sqrt{\frac{nh}{2f(z)\log(1/h)}}\ ||\ g_{n,h,z} - g\ ||_{[0,1]^d} = 0, \tag{9}$$

where $\| g \|_{[0,1]^d} := \sup\{| g(s) |,\ s \in [0,1]^d\}$, and where $\mathcal{S}$ stands for the following Strassen type set of functions mapping $[0,1]^d$ to $\mathbb{R}$:

$$\mathcal{S} := \left\{ g : s \mapsto \int_{[s,1]} \dot{g}(u)du, \text{ for some function } \dot{g} \text{ fulfilling } \int_{[0,1]^d} \dot{g}(u)^2 du. \right\}. \tag{10}$$

Denote by $f_n(K,h,z)$ the Parzen-Rosenblatt density estimator, namely

$$f_n(K,h,z) = \frac{1}{nh} \sum_{i=1}^{n} K\left( \frac{Z_i - z}{h^{1/d}} \right).$$

The main interest of deriving (8) and (9) from Theorem 1 is that it enables us to straightforwardly derive asymptotic confidence bands for $f_n(K,h,z)$ that are uniform in $h \in [h_n, \mathfrak{h}_n]$, which is the subject of the following corollary.

**Corollary 1.1.** *Let $K$ be a kernel with compact support and bounded variation. If both $(h_n)_{n\geq 1}$ and $(\mathfrak{h}_n)_{n\geq 1}$ do satisfy assumptions (HV1)–(HV3), then we have almost surely*

$$\lim_{n\to\infty} \sup_{h\in[h_n,\mathfrak{h}_n]} \sup_{z\in H} \frac{\sqrt{nh}\Big(f_n(K,h,z) - \mathbb{E}\Big(f_n(K,h,z)\Big)\Big)}{\sqrt{2\log(1/h) f(z)}} = \sqrt{\int_{\mathbb{R}^d} K^2 d\lambda},$$

$$\lim_{n\to\infty} \inf_{h\in[h_n,\mathfrak{h}_n]} \inf_{z\in H} \frac{\sqrt{nh}\Big(f_n(K,h,z) - \mathbb{E}\Big(f_n(K,z,h)\Big)\Big)}{\sqrt{2\log(1/h) f(z)}} = -\sqrt{\int_{\mathbb{R}^d} K^2 d\lambda}.$$

*Proof.* By a change of scale, we can assume that $K$ has his support included in $[0,1]^d$. Define the following application that maps the space of all bounded real functions on $[0,1]^d$ to $\mathbb{R}$:

$$\mathcal{R} : g \mapsto \int_{[0,1]^d} g(s) dK(s).$$

Obviously $\mathcal{R}$ is continuous with respect to $\| \cdot \|_{[0,1]^d}$, since $K$ has bounded variation. Noticing that, by an integration by parts, we have

$$f_n(K,h,z) - \mathbb{E}\big(f_n(K,h,z)\big) = \int_{[0,1]^d} g_n(s) dK(s), \tag{11}$$

we readily infer the claimed result, by optimising $\mathcal{R}$ on the limit set $\mathcal{S}$. We omit details. □

In order to state our next corollary, we need to introduce some more notations. Given a positive random variable $h_n^*$, we define

$$\widetilde{\mathbb{E}}\Big(f_n(K,z,h_n^*)\Big) := \int_{\mathbb{R}^d} K\left( \frac{y - z}{{h_n^*}^{1/d}} \right) d\mathbb{P}^{Z_1}(u). \tag{12}$$

**Corollary 1.2.** *Let $h_n^*$ be a sequence of positive random variables satisfying, with probability one:*

$$0 < \liminf_{n\to\infty} \frac{\log(1/h_n)}{\log n} \leq \limsup_{n\to\infty} \frac{\log(1/h_n)}{\log n} < 1. \tag{13}$$

*Then we have almost surely*

$$\lim_{n\to\infty} \sup_{z\in H} \frac{\sqrt{nh_n^*}\Big(f_n(K,z,h_n^*) - \widetilde{\mathbb{E}}\Big(f_n(K,z,h_n^*)\Big)\Big)}{\sqrt{2\log(1/h_n^*)f(z)}} = \sqrt{\int_{\mathbb{R}^d} K^2 d\lambda}, \tag{14}$$

$$\lim_{n\to\infty} \inf_{z\in H} \frac{\sqrt{nh_n^*}\Big(f_n(K,z,h_n^*) - \widetilde{\mathbb{E}}\Big(f_n(K,z,h_n^*)\Big)\Big)}{\sqrt{2\log(1/h_n^*)f(z)}} = -\sqrt{\int_{\mathbb{R}^d} K^2 d\lambda}. \tag{15}$$

*Proof.* The proof is a direct consequence of corollary 1, by manipulating the following countable collection of events

$$F_{r,r'} := \{n^{-1+r} < h_n^* < n^{-1+r'}, \text{ for all large } n\},\ r, r' \in \mathbb{Q} \cap [0,1].$$

We omit details. □

*Remark* 1. Corollaries 1.1 and 1.2 should be compared to Theorem 1 of Einmahl and Mason (2005). Our results rely on assumptions that are stronger than those stated in Einmahl and Mason (2005). However, we derive uniform rates of convergence that are exact and explicit.

*Remark* 2. Several data-driven bandwidths have a well-known asymptotic limit behavior. Typical examples are the bandwidth selectors of Park and Marron (see Park and Marron (1985)) and of Sheather and Jones (see Sheather and Jones (1991)). But the *almost sure* limit behavior of $h_n^*$ is seldom known. However, limits *in probability* are often provided in the literature. For example, it has been proved (see Park and Marron (1985)) that the bandwidth selector $h_{n,1}^*$ of Park and Marron satisfies, under mild conditions,

$$n^{4/13}\left(\frac{h_{n,1}^*}{h_{n,0}} - 1\right) = O_{\mathbb{P}}(1).$$

Here $O_{\mathbb{P}}(1)$ means that the sequence is bounded in probability, and $h_{n,0}$ is the (deterministic) minimizer of

$$\mathbb{E}\Big(\int_{\mathbb{R}} \big(f_n(x) - f(x)\big)^2 dx\Big),$$

which in turn is equivalent to $Cn^{-1/5}$ for some constant $C$. Despite such type of asymptotic results for $h_n^*$ do not meet the requirements of corollary 1.2, it is possible to adapt the latter corollary to derive weaker versions of (14) and (15), which hold *in probability* instead of *almost surely.*

## 2. Proof of Theorem 1

### 2.1. The main tools

Our proof of Theorem 1 relies on two crucial tools. First we shall make use of a criterion in large deviation theory for functional spaces. This criterion, which is mostly due to Arcones (2003), is stated in § 2.1.1. We shall also make use of a concentration inequality (se Proposition 2.2) which is proved by borrowing the arguments of Einmahl and Mason (2000).

#### 2.1.1. Uniform large deviation principles

In the proof of Theorem 1, we shall require large deviation results that are uniform in the rows for triangular arrays of processes. This required uniformity leads us to first state a result that can be straightforwardly derived from Theorem 3.1 of Arcones (2003). In the sequel, $(\epsilon_{n,i})_{n\geq 1, i\leq p_n}$ will always denote a triangular array of positive numbers satisfying $\lim_{n\to\infty} \max_{i\leq p_n} \epsilon_{n,i} = 0$. Given a set $T$, $\mathcal{B}(T)$ will denote the space of bounded real functions on $T$. We shall endow this space with the usual sup-norm $||\cdot||_T$. Let $(E,\vartheta)$ be a topological space. A real function $J : E \to \mathbb{R}^+$ is said to be a rate function (implicitly for $(E,\vartheta)$) when the sets of the kind $\{x \in E : \ J(x) \leq a\}, \ a \geq 0$, are compacts sets of $(E,\vartheta)$. Finally, let $\big(X_{n,i}\big)_{n\geq 1, i\leq p_n}$ be a triangular array of random elements (not necessarily Borel) taking value in $E$. We say that $\big(X_{n,i}\big)_{n\geq 1, i\leq p_n}$ satisfies the uniform large deviation principle (ULDP) for the triangular array $\big(\epsilon_{n,i}\big)_{n\geq 1, i\leq p_n}$ and the rate function $J$ whenever

- For any open set $O \in \mathcal{T}$ we have
$$\liminf_{n\to\infty} \min_{i\leq p_n} \epsilon_{n,i} \log\Big(\mathbb{P}_*\big(X_{n,i} \in O\big)\Big) \geq -J(O).$$
- For any closed set $F \in \mathcal{T}$ we have
$$\limsup_{n\to\infty} \max_{i\leq p_n} \epsilon_{n,i} \log\Big(\mathbb{P}^*\big(X_{n,i}(\cdot) \in F\big)\Big) \leq -J(F).$$

**Remark 2.1.1.**

When referring to outer and inner probabilities $\mathbb{P}^*$ and $\mathbb{P}_*$, one should take care of the underlying probability space, which is taken to be the canonical probability space in our context.

By assumption $(HK5)$ we shall only manipulate true probabilities in our proof of Theorem 1. The following result, which can be seen as a direct corollary of Theorem 2.1 of Arcones (2003), will be used when establishing Proposition 2.3 in the sequel.

**Proposition 2.1.** *Let $(X_{n,i})_{n\geq 1,\ i\leq p_n}$ be a triangular array of random elements of $\mathcal{B}(T)$, and let $\big(\epsilon_{n,i}\big)_{n\geq 1,\ i\leq p_n}$ be a triangular array of positive numbers. Suppose that the following conditions are satisfied.*

1. *There exists a semi distance $\rho$ on $T$ that makes $T$ totally bounded.*
2. *For any $p \geq 1$, and $(t_1, \ldots, t_p) \in \mathbb{R}^p$, the triangular array of random variables $\big(X_{n,i}(t_1), \ldots, X_{n,i}(t_p)\big)_{n\geq 1,\ i\leq p_n}$ satisfies the ULDP for $(\epsilon_{n,i})_{n\geq 1,\ i\leq p_n}$ and a rate function $J_{t_1,\ldots,t_p}$ on $\mathbb{R}^p$.*
3. *For any $\alpha > 0$ and M¿0, there exists $\eta > 0$ such that*

$$\limsup_{n\to\infty} \max_{i\leq p_n} \epsilon_{n,i} \log\Big(\mathbb{P}^*\Big(\sup_{t,s:\ \rho(t,s)\leq\eta} |X_{n,i}(t) - X_{n,i}(s)| > \alpha\Big)\Big) \leq -M. \tag{16}$$

*Then $(X_{n,i}(\cdot))_{n\geq 1,\ i\leq p_n}$ satisfies the ULDP in $\big(\mathcal{B}(T), ||\cdot||_T\ \big)$ for $(\epsilon_{n,i})_{n\geq 1,\ i\leq p_n}$ and the following rate function*

$$J\big(\Psi\big) := \sup_{p\geq 1,\ (t_1,\ldots,t_p)\in\mathbb{R}^p} \Big\{J_{t_1,\ldots,t_p}\big(\Psi(t_1), \ldots, \Psi(t_p)\big)\Big\},\ \Psi \in \mathcal{B}(T).$$

*Proof.* The proof is a direct copy of the proof of Theorem 2.1 in Arcones (2003), replacing $\mathbb{P}^*\big(U_n \in F\big) \leq \ldots$ by $\max_{i\leq p_n} \mathbb{P}^*\big(X_{n,i} \in F\big) \leq \ldots$ and $\mathbb{P}_*\big(U_n \in O\big) \geq \ldots$ by $\min_{i\leq p_n} \mathbb{P}_*\big(X_{n,i} \in O\big) \geq \ldots$, and so on. We avoid writing the proof for this reason. □

*2.1.2. A concentration inequality*

For any real Borel function $g$ we set

$$T_n(g) = \sum_{i=1}^n g(Z_i) - \mathbb{E}\big(g(Z_i)\big),\ g \in \mathcal{F}. \tag{17}$$

The following concentration inequality for local empirical processes will play a crucial role in the sequel. It states a somewhat finite distance version of the arguments of Einmahl and Mason (2005).

**Proposition 2.2.** *Let $\mathcal{F}$ be a class of functions on $\mathbb{R}^d$ with measurable envelope function $F$ satisfying, for some constants $\tau > 0$ and $h \in (0,1)$,*

$$\sup_{g\in\mathcal{F}} \mathrm{Var}\big(g(Z_1)\big) \leq \tau^2 h.$$

*Assume that there exists $C, v, \beta_0 > 0$ and $p > 2$ fulfilling, for all $0 < \epsilon < 1$,*

$$\begin{aligned}\mathcal{N}(\epsilon, \mathcal{F}) &\leq C\epsilon^{-v},\\ \mathbb{E}\big(F(Y)^2\big) &\leq \beta_0^2.\end{aligned}$$

*Then there exists a universal constant $A_2 > 0$ and a parameter $D_1(v) > 0$ depending only on $v$ such that, for fixed $\rho_0 > 0$, $p > 2$ and $\delta_0$, if $h > 0$ satisfies*

$$C_1 := \max\Big\{1, \big(4\delta_0\sqrt{v+1}/\tau\big)^{\frac{1}{1/2-1/p}}, \big(\rho_0\delta_0/\tau^2\big)^{\frac{1}{1/2-1/p}}\Big\} \leq \frac{nh}{\log(1/h)}, \tag{18}$$

$$C_2 := \min\{1/\tau^2\beta_0, \tau^2\} \geq h \text{ and} \tag{19}$$

$$\sup_{\substack{g\in\mathcal{F},\\ z\in\mathbb{R}^d}} \mid g(z) \mid \leq \delta_0(nh/\log(1/h))^{1/p}, \tag{20}$$

*then we have*

$$\mathbb{P}\Big(\max_{1\leq m\leq n} || T_m ||_{\mathcal{F}} \geq (\tau+\rho_0)D_1(nh\log(1/h))^{1/2}\Big)$$
$$\leq 4\exp\bigg(-A_2\Big(\frac{\rho_0}{\tau}\Big)^2\log(1/h)\bigg).$$

*Proof.* By (18) and (20) we have

$$\begin{aligned}\sup_{g\in\mathcal{G},\ z\in\mathbb{R}^d} \mid g(z) \mid \leq& \delta_0\Big(\frac{nh}{\log(1/h)}\Big)^{1/p}\\ \leq& \frac{1}{2\sqrt{v+1}}\sqrt{\frac{n\tau^2 h}{2\log(1/h)}}\\ \leq& \frac{1}{2\sqrt{v+1}}\sqrt{\frac{n\tau^2 h}{\log(\beta_0\vee 1/\tau^2 h)}}.\end{aligned} \tag{21}$$

Here, (21) is a consequence of (19).

Denote by $(\epsilon_i)_{i\geq 1}$ an i.i.d. sequence of random variables independent of $(Z_i)_{i\geq 1}$ with $\mathbb{P}(\epsilon_1 = \pm 1) = 1/2$. Applying Proposition A.2 of Einmahl and Mason (2000) with $\beta := \beta_0^2$ and $\sigma^2 := \tau^2 h$ we get, for a universal constant $A_3 > 0$,

$$\begin{aligned}\mu_n(\mathcal{F}) :=& \mathbb{E}\bigg(\sup_{g\in\mathcal{F}}\Big|\sum_{i=1}^n \epsilon_i g(Z_i)\Big|\bigg)\\ \leq& A_3\sqrt{vn\tau^2 h\log(1/\tau^2 h)}\\ \leq& \sqrt{2v}A_3\tau\sqrt{nh\log(1/h)}.\end{aligned} \tag{22}$$

We shall now apply inequality A.1 in Einmahl and Mason (2005). According to the notations of that inequality, we choose $M := \delta_0(nh/\log(1/h))^{1/p}$. We then have, writing $D_1 := \max\{A_1, A_1A_3\sqrt{2v}\}$,

$$\mathbb{P}\Big(\max_{1\leq m\leq n} || T_m ||_{\mathcal{F}} \geq (\tau+\rho_0)D_1(nh\log(1/h))^{1/2}\Big)$$
$$\leq \mathbb{P}\bigg(\max_{1\leq m\leq n} || T_m ||_{\mathcal{F}} \geq A_1\Big(\rho_0(nh\log(1/h))^{1/2} + \mu_n(\mathcal{F})\Big)\bigg) \tag{23}$$
$$\leq 2\bigg[\exp\Big(-\frac{A_2\rho_0^2 nh\log(1/h)}{n\tau^2 h}\Big) + \exp\Big(-\frac{A_2\rho_0(nh\log(1/h))^{1/2}}{\delta_0(nh/\log(1/h))^{1/p}}\Big)\bigg] \tag{24}$$
$$\leq 4\exp\Big(-\frac{A_2\rho_0^2}{\tau^2}\log(1/h)\Big). \tag{25}$$

Here, (24) is a direct application of inequality A.1 in Einmahl and Mason (2005), while inequality (25) is a consequence of (18). This concludes the proof of Proposition 2.2. □

### 2.2. Proof of part (i) of Theorem 1

We will make repeatedly use of the following obvious argument:

$$\sup_{z\in H} f(z)^{-1/2} =: \beta < \infty. \tag{26}$$

First select $\epsilon > 0$ arbitrarily. We claim that, almost surely,

$$\limsup_{n\to\infty} \sup_{\substack{h_n\le h\le \mathfrak{h}_n,\\ z\in H}} \inf_{\Psi\in\mathbb{K}} \left|\left| \frac{G_n(\cdot,h,z)}{\sqrt{2f(z)nh\log(1/h)}} - \Psi \right|\right|_{\mathcal{G}} \le \epsilon. \tag{27}$$

To prove this, we introduce some parameters that will be properly adjusted in the sequel. Recall that $\alpha > 0$ appears in (4). Let $\gamma > 0$, $\rho > 1$ and $0 < \delta < \alpha/4$ be real numbers. We shall invoke some usual blocking arguments along the subsequence $n_k := [(1+\gamma)^k]$, $k \ge 1$. For fixed $k \ge 1$, consider the following discretisation of $[h_{n_k}, \mathfrak{h}_{n_{k-1}}]$.

$$h_{n_k,R_k} := \mathfrak{h}_{n_{k-1}},\ \ h_{n_k,l} := \rho^l h_{n_k},\ l = 0,\dots,R_k - 1, \tag{28}$$

where $R_k := [\log(\mathfrak{h}_{n_{k-1}}/h_{n_k})/\log(\rho)] + 1$, and $[u]$ denotes the only integer $q$ fulfilling $q \le u < q+1$. Since $h_n$ and $\mathfrak{h}_n$ satisfy (HV1)–(HV3), the triangular array $h_{n_k,l}$, $0 \le l \le R_k$ satisfies the two following properties that will play a crucial role in our arguments (see our *key argument 1* below).

$$\lim_{k\to\infty} \max_{0\le l\le R_k} h_{n_k,l} = 0 \tag{29}$$

$$\lim_{k\to\infty} \min_{0\le l\le R_k} \frac{n_k h_{n_k,l}}{\log(1/h_{n_k,l})} = \infty. \tag{30}$$

Recall that $I^d := [0,1]^d$. For each $0 \le l \le R_k$, we proceed as in Mason (2004), covering $H$ by pairwise disjoint hypercubes written as

$$\Gamma_{k,l,j} := \Big\{ z_{k,l,j} + \big[0, (\delta h_{n_k,l})^{1/d}\big)^d \Big\},\ 1 \le j \le J_l, \tag{31}$$

with $z_{k,l,j} \in H$. Since $0 < \delta < \alpha/4$, we have, for all $k \ge 1$ and $0 \le l \le R_k$,

$$H \subset \bigcup_{j=1}^{J_l} \Gamma_{k,l,j} \subset H^{\alpha/2}. \tag{32}$$

Notice that by construction:

$$J_l \le \frac{C}{h_{n_k,l}},\ k \ge 1,\ 0 \le l \le R_k, \tag{33}$$

where $C := C(\delta)$ depends on $\delta > 0$ and on the volume of $H$ only. Set $N_k := \{n_{k-1}+1,\dots,n_k\}$ whenever $n_{k-1} < n_k$, and $N_k = \emptyset$ elsewhere. For $A \subset \mathcal{B}(\mathcal{G})$ and $\epsilon > 0$ we write the $\epsilon$ neighbourhood of $A$ as

$$A^\epsilon := \Big\{ \Psi \in \mathcal{B}(\mathcal{G}),\ \inf_{\Psi'\in A} \|\ \Psi - \Psi'\ \|_{\mathcal{G}} < \epsilon \Big\}. \tag{34}$$

For all large $k$, we split the following probabilities in two.

$$\begin{aligned}
\mathbb{P}_k &:= \mathbb{P}\Bigg( \bigcup_{\substack{n\in N_k, z\in H,\\ h\in[h_n,\mathfrak{h}_n]}} \Bigg\{ \inf_{\Psi\in\mathbb{K}} \Bigg\| \frac{G_n(\cdot,h,z)}{\big(2f(z)nh\log(1/h)\big)^{1/2}} - \Psi \Bigg\|_{\mathcal{G}} > 4\epsilon \Bigg\} \Bigg) \\
&\leq \mathbb{P}\Bigg( \bigcup_{\substack{n\in N_k, z\in H,\\ h\in[h_{n_k},\mathfrak{h}_{n_{k-1}}]}} \Bigg\{ \frac{G_n(\cdot,h,z)}{\big(2f(z)nh\log(1/h)\big)^{1/2}} \notin \mathbb{K}^{4\epsilon} \Bigg\} \Bigg) \\
&\leq \mathbb{P}\Bigg( \bigcup_{\substack{n\in N_k,\ 0\leq l\leq R_k,\\ 1\leq j\leq J_l}} \Bigg\{ \frac{G_n(\cdot,h_{n_k,l},z_{k,l,j})}{\big(2f(z_{k,l,j})n_k h_{n_k,l}\log(1/h_{n_k,l})\big)^{1/2}} \notin \mathbb{K}^{2\epsilon} \Bigg\} \Bigg) \\
&\quad + \mathbb{P}\Bigg( \max_{\substack{n\in N_k, 0\leq l\leq R_k-1,\\ 1\leq j\leq J_l}} \sup_{\substack{h_{n_k,l}\leq h\leq \rho h_{n_k,l},\\ z\in\Gamma_{k,l,j}}} \Bigg\| \frac{G_n(\cdot,h,z)}{\big(2f(z)nh\log(1/h)\big)^{1/2}} \\
&\qquad - \frac{G_n(\cdot,h_{n_k,l},z_{k,l,j})}{\big(2f(z_{k,l,j})n_k h_{n_k,l}\log(1/h_{n_k,l})\big)^{1/2}} \Bigg\|_{\mathcal{G}} > 2\epsilon \Bigg) \\
&=: \mathbb{P}_{1,k} + \mathbb{P}_{2,k}. \qquad (35)
\end{aligned}$$

Our aim is to prove that $\mathbb{P}_{1,k}$ and $\mathbb{P}_{2,k}$ are both summable in $k$, which would prove part (i) of Theorem 1 by an application of Borel-Cantelli's lemma to $\mathbb{P}_k$.

#### 2.2.1. Step 1: blocking and poissonisation

By a blocking argument that is similar to Ottaviani's inequality (see for example Deheuvels and Mason (1992), Lemma 3.4), we have, for fixed $k \geq 1$, $l \leq R_k$, $j \leq J_l$:

$$\begin{aligned}
&\mathbb{P}\Bigg( \bigcup_{n\in N_k} \Bigg\{ \frac{G_n(\cdot,h_{n_k,l},z_{k,l,j})}{(2f(z_{k,l,j})n_k h_{n_k,l}\log(1/h_{n_k,l})))^{1/2}} \notin \mathbb{K}^{2\epsilon} \Bigg\} \Bigg) \\
&\qquad \leq \frac{1}{m_k}\mathbb{P}\Bigg( \frac{G_{n_k}(\cdot,h_{n_k,l},z_{k,l,j})}{(2f(z_{k,l,j})n_k h_{n_k,l}\log(1/h_{n_k,l}))^{1/2}} \notin \mathbb{K}^{\epsilon} \Bigg), \qquad (36)
\end{aligned}$$

where

$$m_k := \min_{\substack{l\leq R_k,\ j\leq J_l,\\ 0\leq m\leq n_k-n_{k-1}}} \mathbb{P}\Bigg( \Bigg\| \frac{G_m(\cdot,h_{n_k,l},z_{k,l,j})}{\big(2f(z_{k,l,j})n_k h_{n_k,l}\log(1/h_{n_k,l})\big)^{1/2}} \Bigg\|_{\mathcal{G}} \leq \epsilon \Bigg).$$

To control $m_k$, we shall invoke an argument that will be repeatedly used in that article. Roughly speaking, we make use of the arguments of Mason (2004) replacing $h_{n_k}$ by $h_{n_k,l}$. This leads us to consider the following classes of function, for $k \geq 1, l \leq R_k, j \leq J_l$:

$$\mathcal{F}_{k,l,j} := \Bigg\{ f(z_{i,n_k})^{-1/2} K\Bigg( \frac{\cdot - z_{k,l,j}}{h_{n_k,l}{}^{1/d}} \Bigg),\ K\in\mathcal{G} \Bigg\} \subset \mathcal{F}.$$

To obtain an upper bound that holds uniformly in $h_{n_k,l}$, we shall show that all these classes do satisfy the assumptions of Proposition 2.2 simultaneously in $h_{n_k,l}$ an $z_{k,l,j}$. To prove this, assertions (29) and (30) will play a crucial role.

***Key argument 1***: By Lemma 1 in Mason (2004) (Bochner's lemma), and by (29) we have

$$\sup_{g\in\mathcal{F}_{k,l,j}} \operatorname{Var}\big(g(Z_1)\big) \leq h_{n_k,l}(1+v_k),$$

where $v_k \to 0$ as $k\to\infty$. Notice that each $\mathcal{F}_{k,l,j}$ is uniformly bounded by 1, in virtue of (HK2) and (26). We now use the notations of Proposition 2.2 with $\tau := 2$, $\rho_0 = 2$, $\delta_0 = 1$, $p = 4$, $n := n_k - n_{k-1}$ and the constants $C, v$ appearing in assumption (HK4). Since $n_k - n_{k-1} \sim \gamma/(1+\gamma)n_k$ and by both (30) and (29), we have, for all large $k$, $h_{n_k,l} \leq C_2$ and

$$\min_{0\leq l\leq R_k} \left(\frac{(n_k-n_{k-1})h_{n_k,l}}{\log(1/h_{n_k,l})}\right)^{1/2} \geq C_1,$$

$$\max_{\substack{0\leq l\leq R_k,\\ 1\leq j\leq J_l}} \sup_{\substack{g\in\mathcal{F}_{k,l,j},\\ z\in\mathbb{R}^d}} \mid g(z) \mid \leq \min_{0\leq l\leq R_k} 2\Big((n_k-n_{k-1})h_{n_k,l}/\log(1/h_{n_k,l})\Big)^{1/p}. \quad (37)$$

This implies that, for all large $k$, *all the classes* $\mathcal{F}_{k,l,j}$ do satisfy the assumptions of Proposition 2.2. Moreover, for $\gamma > 0$ small enough and for $k$ large enough we have

$$4D_1(v)(n_k - n_{k-1})^{1/2} \leq \frac{\epsilon}{\beta} n_k^{1/2},$$

and hence

$$\mathbb{P}\Big(\max_{m\leq n_k-n_{k-1}} \|\, G_m(\cdot, h_{n_k,l}, z_{k,l,j})\,\|_{\mathcal{G}} \geq \epsilon\big(2f(z_{k,l,j})n_k\log(1/h_{n_k,l})\big)^{1/2}\Big)$$

$$\leq \mathbb{P}\Big(\max_{m\leq n_k-n_{k-1}} \|\, T_m\,\|_{\mathcal{F}} \geq (\rho+\tau)D_1(v)\big(2(n_k-n_{k-1})h_{n_k,l}\log(1/h_{n_k,l})\big)^{1/2}\Big)$$

$$\leq 4\exp\big(-A_2\log(1/\mathfrak{h}_{n_{k-1}})\big).$$

Therefore, for $\gamma > 0$ small enough and for $k$ large enough we have $m_k \geq 1/2$. Now let

$$\widetilde{G}_n(K,h,z) := \sum_{i=1}^{\eta_n} K\left(\frac{Z_i - z}{h^{1/d}}\right) - \mathbb{E}\left(\frac{Z_i - z}{h^{1/d}}\right) \quad (38)$$

be the Poissonized version of $G_n$. Here, $(\eta_n)_{n\geq 1}$ is a Poisson random variable with expectation $n$, and independent of $(Z_i)_{i\geq 1}$. Recalling that by construction, there exists $C = C(\delta) < \infty$ such that $J_l \leq C(\delta)/h_{n_k,l}$, $l = 1, \ldots, R_k$, it follows that, ultimately as $k\to\infty$,

$$\mathbb{P}_{1,k} \leq 2\sum_{l=0}^{R_k} \frac{C}{h_{n_k,l}} \max_{j\leq J_l} \mathbb{P}\left(\frac{G_{n_k}(\cdot, h_{n_k,l}, z_{k,l,j})}{\big(2f(z_{k,l,j})n_k h_{n_k,l}\log(1/h_{n_k,l})\big)^{1/2}} \notin \mathbb{K}^{\epsilon}\right)$$

$$\leq 4 \sum_{l=0}^{R_k} \frac{C}{h_{n_k,l}} \max_{j\leq J_l} \mathbb{P}\left( \frac{\widetilde{G}_{n_k}(\cdot, h_{n_k,l}, z_{k,l,j})}{\left(2f(z_{k,l,j})n_k h_{n_k,l}\log(1/h_{n_k,l})\right)^{1/2}} \notin \mathbb{K}^{\epsilon} \right)$$

$$=: 4 \sum_{l=0}^{R_k} \frac{C}{h_{n_k,l}} \max_{j\leq J_l} \mathbb{P}_{k,l,j}. \tag{39}$$

The last inequality is a consequence of usual poissonization inequalities (see, e.g., Giní et al. (2003), Lemma 2.1).

#### *2.2.2. Step 2: A uniform large deviation result*

In order to control the $\mathbb{P}_{k,l,j}$uniformly in $l$ and $j$, we shall establish a uniform large deviation principle that is stated in the next proposition. Recall that $J$ has been defined in (3). Some routine analysis shows that $J$ is a rate function on $\mathcal{B}(\mathcal{G})$.

**Proposition 2.3.** *Let $(n_k)_{k\geq 1}$ be a strictly increasing integer-valued sequence and $\left(z_{k,l,j}\right)_{k\geq 1, l\leq R_k,\, j\leq J_l}$ a triangular array of points belonging to $H$. Under (HV1)–(HV3) and (HK1)–(HK5), the triangular array*

$$\left( \frac{\widetilde{G}_{n_k}(\cdot, h_{n_k,l}, z_{k,l,j})}{\left(2f(z_{k,l,j})n_k h_{n_k,l}\log(1/h_{n_k,l})\right)^{1/2}} \right)_{k\geq 1, l\leq R_k, j\leq J_l},$$

*satisfies the ULDP for the rate function $J$ and the triangular array*

$$\epsilon_{k,l,j} := \left(\log(1/h_{n_k,l})\right)^{-1},\ k\geq 1, l\leq R_k, j\leq J_l.$$

*Proof.* To prove Proposition 2.3, we shall make use of Proposition 2.1, and we hence have to check conditions 1, 2 and 3 of that theorem. Compared to Proposition 1 of Mason (2004), the present proposition adds a uniformity in $h_{n_k,l}$. Checking condition 2 of Proposition 2.1 readily follows the lines of the proof of Mason according to the following remarks: We can apply his Fact 2 with an additional uniformity in $h_{n_k,l}$, as $\mathfrak{h}_n \to 0$ and hence Bochner's lemma still holds uniformly in $h_{n_k,l}$. We can also apply his Fact 3, replacing $h_n$ by $h_{n_k,l}$. Hence, his assertion (4.16) still holds replacing $h_n$ by $h_{n_k,l}$, with an additional uniformity in $h_{n_k,l}$. Now define the following distance on $\mathcal{G}$.

$$d^2(K, K') := \int_{\mathbb{R}^d} \left(K - K'\right)^2(z)dz.$$

It remains to show that for any $M > 0, \overline{\alpha} > 0$, there exists $\overline{\delta} > 0$ fulfilling

$$\limsup_{k\to\infty} \max_{\substack{0\leq l\leq R_k,\\ 1\leq j\leq J_l}} \epsilon_{k,l,j} \times \log\left( \mathbb{P}\left( \sup_{\substack{K,K'\in\mathcal{G},\\ d(K,K')\leq\overline{\delta}}} \left| \frac{\widetilde{G}_{n_k}(K, h_{n_k,l}, z_{k,l,j}) - \widetilde{G}_{n_k}(K', h_{n_k,l}, z_{k,l,j})}{\sqrt{2f(z_{k,l,j})h_{n_k,l}\log(1/h_{n_k,l})}} \right| \geq \overline{\alpha} \right) \right) \leq -M, \tag{40}$$

So as conditions 1 and 3 of Proposition 2.1 would be checked. Choose $M > 0$ and $\overline{\alpha} > 0$ arbitrarily. For each $k \geq 1, 0 \leq l \leq R_k, 1 \leq j \leq J_l$ and $\overline{\delta} > 0$, consider the following class of functions:

$$\mathcal{F}_{k,j,l,\overline{\delta}} := \left\{ f(z_{k,l,j})^{-1/2}(K - K')\left(\frac{\cdot - z_{k,l,j}}{h_{n_k,l}{}^{1/d}}\right),\ d^2(K, K') \leq \overline{\delta} \right\}.$$

Let $D_1(2v)$ be the as in Proposition 2.2 (recall that $v > 0$ appears in assumption (HK4)). We have, for any $k \geq 1,\ 0 \leq l \leq R_k,\ 1 \leq j \leq J_l$,

$$\begin{aligned}
\widetilde{\mathbb{P}}_{k,l,j} &:= \mathbb{P}\Bigg( \sup_{\substack{K,K' \in \mathcal{G}, \\ d(K,K') \leq \overline{\delta}}} \left| \frac{\widetilde{G}_{n_k}(K, h_{n_k,l}, z_{k,l,j}) - \widetilde{G}_{n_k}(K', h_{n_k,l}, z_{k,l,j})}{\sqrt{2 f(z_{k,l,j}) h_{n_k,l} \log(1/h_{n_k,l})}} \right| \geq 2D_1(2v)\overline{\alpha} \Bigg) \\
&= \mathbb{P}\Bigg( \left|\left| \frac{\widetilde{G}_{n_k}(\cdot, h_{n_k,l}, z_{k,l,j})}{\big(2n_k h_{n_k,l} \log(1/h_{n_k,l})\big)^{1/2}} \right|\right|_{\mathcal{F}_{k,l,j,\overline{\delta}}} \geq 2D_1(2v)\overline{\alpha} \Bigg) \\
&= \sum_{m=1}^{\infty} \mathbb{P}(\eta_{n_k} = m)\mathbb{P}\Bigg( \left|\left| \frac{G_m(\cdot, h_{n_k,l}, z_{k,l,j})}{\big(2n_k h_{n_k,l} \log(1/h_{n_k,l})\big)^{1/2}} \right|\right|_{\mathcal{F}_{k,l,j,\overline{\delta}}} \geq 2D_1(2v)\overline{\alpha} \Bigg) \\
&\leq \mathbb{P}\Bigg( \max_{m=1,\dots,2n_k} \left|\left| \frac{G_m(\cdot, h_{n_k,l}, z_{k,l,j})}{\big(2n_k h_{n_k,l} \log(1/h_{n_k,l}))\big)^{1/2}} \right|\right|_{\mathcal{F}_{k,l,j,\overline{\delta}}} \geq 2D_1(2v)\overline{\alpha} \Bigg) \\
&\quad + \mathbb{P}\big(\eta_{n_k} > 2n_k\big) \\
&=: \widetilde{\mathbb{P}}_{k,l,j,1} + \widetilde{\mathbb{P}}_{k,l,j,2}.
\end{aligned} \tag{41}$$

By Chernoff's inequality we have

$$\widetilde{\mathbb{P}}_{k,l,j,2} \leq \exp\big(-(2\log 2 - 1)n_k\big),\ 0 \leq l \leq R_k,\ 1 \leq j \leq J_l. \tag{42}$$

By (HK4) we have, by simple arguments,

$$\mathcal{N}(\epsilon, \mathcal{F}_{k,j,l,\overline{\delta}}) \leq \big(\mathcal{N}(\epsilon/2\beta, \mathcal{F})\big)^2 \leq \big(2\beta\big)^{2v} C^2 \epsilon^{-2v} =: C'\epsilon^{-v'},\ 0 < \epsilon < 1. \tag{43}$$

An application of Lemma 1 in Mason (2004) in combination with (29) leads to the following inequality, for all large $k$:

$$\max_{\substack{0 \leq l \leq R_k, \\ 1 \leq j \leq J_l}} \sup_{g \in \mathcal{F}_{k,j,l,\overline{\delta}}} h_{n_k,l}^{-1} \mathrm{Var}\left( g\left( \frac{Z_1 - z_{k,l,j}}{h_{n_k,l}^{1/d}} \right) \right) \leq 2\overline{\delta} \tag{44}$$

Reasoning as in the *key argument 1*, we conclude that, for all large $k$, each class $\mathcal{F}_{k,l,j,\overline{\delta}},\ 0 \leq l \leq R_k,\ 1 \leq j \leq J_l$ do fulfill conditions (20), (18) and (19) of Proposition 2.2 with $\rho_0 := \overline{\alpha},\ \tau := \sqrt{2\overline{\delta}},\ \delta_0 := 2,\ p := 4,\ n := 2n_k$ and $C', v'$ appearing in (43). Applying Proposition 2.2, we get, for all large $k$ and for each $0 \leq l \leq R_k,\ 1 \leq j \leq J_l$,

$$\begin{aligned}&\widetilde{\mathbb{P}}_{k,l,j,1}\\ &\leq \mathbb{P}\Big(\max_{m=1,\dots,2n_k} ||\, T_m\, ||_{\mathcal{F}_{k,l,j,\overline{\delta}}} \geq D_1(v')\big(\overline{\alpha}+\sqrt{2\overline{\delta}}\big)\sqrt{2n_k f(z_{k,l,j})h_{n_k,l}\log(1/h_{n_k,l})}\Big)\\ &\leq 4\exp\Big(-A_2\frac{\overline{\alpha}^2}{2\overline{\delta}}\log(1/h_{n_k,l})\Big). \qquad (45)\end{aligned}$$

Notice that (45) is true for all $\overline{\delta}>0$ satisfying $2\overline{\delta}\leq\overline{\alpha}$. By (42) in conjunction with (45), we have for $\overline{\delta}>0$ small enough, ultimately as $k\to\infty$,

$$\begin{aligned}\widetilde{\mathbb{P}}_{k,l,j} \leq& 4\exp\big(-M\log(1/h_{n_k,l})\big)+\exp\big(-(2\log 2-1)n_k\big)\\ \leq& 5\exp\big(-M\log(1/h_{n_k,l})\big).\end{aligned}$$

This shows that (40) is true. We now refer to Arcones (200), Theorem 4.2, for the proof of the fact that

$$\sup_{p\geq 1,\ (K_1,\dots,K_p)\in\mathcal{G}^p} J_{K_1,\dots,K_p}(\Psi(K_1),\dots\Psi(K_p)) = J(\Psi),$$

which completes the proof of Proposition 2.3 by an application of Proposition 2.1. □

### *2.2.3. Step 3: summability of* $\mathbb{P}_{1,k}$

For fixed $\epsilon>0,\delta>0,\rho>1$ and for $\gamma>0$ small enough, we apply Proposition 2.3 to the following closed subset of $\big(\mathcal{B}(\mathcal{G}),||\cdot||_{\mathcal{G}}\big)$.

$$F_\epsilon := \Big\{\Psi\in\mathcal{B}(\mathcal{G}),\ \inf_{\Psi'\in\mathbb{K}} ||\ \Psi-\Psi'\ ||_{\mathcal{G}}\geq\epsilon\Big\}.$$

***Key argument 2***: By lower semicontinuity of $J$ on $\big(\mathcal{B}(\mathcal{G}),||\cdot||_{\mathcal{G}}\big)$, there exists $\alpha_1>0$ such that $J\big(F_\epsilon\big)=1+2\alpha_1$. Hence, inequality (39) becomes

$$\begin{aligned}\mathbb{P}_{1,k} \leq& 4\sum_{l=0}^{R_k}\frac{C}{h_{n_k,l}}\exp\Big(-\big(J\big(F_\epsilon\big)-\alpha_1\big)\log(1/h_{n_k,l})\Big)\\ \leq& 4C\sum_{l=0}^{R_k} h_{n_k,l}{}^{\alpha_1} \leq 4Ch_{n_k}^{\alpha_1}\sum_{l=0}^{R_k}\rho^{\alpha_1 l} \leq \frac{4C}{\rho^{\alpha_1}-1}h_{n_k}^{\alpha_1}\rho^{\alpha_1(R_k+1)}.\end{aligned}$$

Recall that $C>0$ depends only on $\delta>0$ and $H$. Since we have by construction $R_k\leq\log(\mathfrak{h}_{n_{k-1}}/h_{n_k})/\log(\rho)+1$, we deduce that

$$\mathbb{P}_{1,k}\leq\frac{4C\rho^{2\alpha_1}}{\rho^{\alpha_1}-1}\mathfrak{h}_{n_{k-1}}^{\alpha_1}.$$

But (HV3) ensures that $\mathfrak{h}_{n_{k-1}}^{\alpha_1}$ is summable, which in turn implies that$(\mathbb{P}_{k,1})_{k\geq 1}$ is widely summable in $k$, whatever the choice of $\delta>0$, $\rho>1$ and $\gamma>0$.

### 2.2.4. Step 4: an upper bound for $\mathbb{P}_{2,k}$

Now for an arbitrary $\epsilon > 0$ we shall adjust $\delta > 0$ and $\rho > 1$ such that, for $\gamma > 0$ small enough, the sequence $\mathbb{P}_{2,k}$ has a finite sum in $k$. We start by the following decomposition.

$$\begin{aligned}
&\mathbb{P}_{2,k} \\
&:= \mathbb{P}\Bigg( \max_{\substack{n\in N_k, 0\le l\le R_k-1,\\ 1\le j\le J_l}} \sup_{\substack{h_{n_k,l}\le h\le \rho h_{n_k,l},\\ z\in\Gamma_{k,l,j}}} \Bigg\| \frac{G_n(\cdot,h,z)}{(2f(z)nh\log(1/h))^{1/2}} \\
&\qquad - \frac{G_n(\cdot,h_{n_k,l},z_{k,l,j})}{(2f(z_{k,l,j})n_k h_{n_k,l}\log(1/h_{n_k,l}))^{1/2}} \Bigg\|_{\mathcal{G}} > 2\epsilon \Bigg) \\
&\le \mathbb{P}\Bigg( \max_{\substack{n\in N_k, 0\le l\le R_k-1,\\ 1\le j\le J_l}} \sup_{\substack{h_{n_k,l}\le h\le \rho h_{n_k,l},\\ z\in\Gamma_{k,l,j}}} \Bigg\| \frac{G_n(\cdot,h,z)-G_n(\cdot,h_{n_k,l},z_{k,l,j})}{\sqrt{2f(z_{k,l,j})n_k h_{n_k,l}\log(1/h_{n_k,l})}} \Bigg\|_{\mathcal{G}} > \epsilon \Bigg) \\
&+ \mathbb{P}\Bigg( \max_{\substack{n\in N_k, l\le R_k-1,\\ 1\le j\le J_l}} \sup_{\substack{h_{n_k,l}\le h\le \rho h_{n_k,l},\\ z\in\Gamma_{k,l,j}}} B_{k,n,h,z} \Bigg\| \frac{G_n(\cdot,h,z)}{\sqrt{2f(z)n_k h_{n_k,l}\log(1/h_{n_k})}} \Bigg\|_{\mathcal{G}} > \epsilon \Bigg) \\
&=: \mathbb{P}_{2,1,k} + \mathbb{P}_{2,2,k},
\end{aligned}$$

Here we have set, for $z \in \Gamma_{k,j,l}$,

$$B_{k,n,h,z} := \left| \sqrt{\frac{2f(z_{k,l,j})n_k h_{n_k,l}\log(1/h_{n_k,l})}{2f(z)nh\log(1/h)}} - 1 \right|. \tag{46}$$

In order to control $\mathbb{P}_{2,2,k}$, we make the following decomposition for all large $k$:

$$\mathbb{P}_{2,2,k} \le \sum_{\substack{0\le l\le R_k-1,\\ 1\le j\le J_l}} \mathbb{P}_{2,2,k,l,j},$$

where

$$\mathbb{P}_{2,2,k,l,j} := \mathbb{P}\Bigg( \max_{n\in N_k} \sup_{\substack{h_{n_k,l}\le h\le \rho h_{n_k,l},\\ z\in\Gamma_{k,l,j}}} B_{k,n,h,z} \Bigg\| \frac{G_n(\cdot,h,z)}{\sqrt{f(z)n_k h\log(1/h_{n_k,l})}} \Bigg\|_{\mathcal{G}} \ge \epsilon \Bigg).$$

Some usual analysis based on (H f), (29), (30), and $n_k - n_{k-1} \sim n_{k-1}\gamma/(1+\gamma)$ show that, for any choice of $\delta > 0$, $\rho > 1$ small enough we have, for all large $k$:

$$\max_{\substack{n\in N_k,\ l\le R_k,\\ j\le J_l}} \sup_{\substack{h_{n_k,l}\le h\le \rho h_{n_k,l},\\ z\in\Gamma_{k,l,j}}} B_{k,n,h,z} \le \gamma(1+2\gamma)+2\gamma. \tag{47}$$

Hence, for any choice of $\gamma > 0$ small enough, we have (recall Proposition 2.2 and assumption (HV4))

$$\liminf_{k\to\infty} \inf_{\substack{n\in N_k,\ z\in H,\\ h_{n_k}\le h\le \mathfrak{h}_{n_{k-1}}}} \epsilon B^{-1}_{k,n,h,z} \ge 4D_1(v). \tag{48}$$

Now consider the following classes of functions:

$$\overline{\mathcal{F}}_{k,l,j} := \left\{ f(z)^{-1/2} K\left(\frac{\cdot - z}{h^{1/d}}\right),\ z \in \Gamma_{k,l,j},\ h_{n_k,l} \le h \le \rho h_{n_k,l} \right\}.$$

By (29) and Lemma 1 in Mason (2004) we have, for all large $k$:

$$\max_{l \le R_k,\ j \le J_l} \sup_{g \in \overline{\mathcal{F}}_{k,l,j}} h_{n_k,l}^{-1} \mathrm{Var}\big(f(Z_1)\big) \le 4. \tag{49}$$

Recall (26). According to (HK4) we have

$$\mathcal{N}\big(\epsilon, \overline{\mathcal{F}}_{k,l,j}\big) \le \mathcal{N}\big(\epsilon/\beta, \mathcal{F}\big) \le C\beta^v \epsilon^{-v} =: C'\epsilon^{-v}.$$

Proceeding similarly a in the *key argument 1*, we infer that all the classes $\overline{\mathcal{F}}_{k,l,j}$, $0 \le l \le R_k$, $1 \le j \le J_l$ do satisfy conditions (20) (18) and (19) in Proposition 2.2, with $\tau := \sqrt{A_2} \wedge 1$, $\rho := 2$, $\delta_0 = 2$, $p = 4$, $n := n_k$, $h := h_{n_k,l}$, $C'$ and $v$. Making use once again of Proposition 2.2, and assuming that $\gamma$ is small enough to fulfill (48) we get that, for all large $k$ and $0 \le l \le R_k$, $1 \le j \le J_l$,

$$\begin{aligned}
\mathbb{P}_{2,2,k,l,j} \le & \mathbb{P}\Big(\max_{n \in N_k} || T_n ||_{\overline{\mathcal{F}}_{k,l,j}} \ge 4D_1(v)\big(2n_k \log(1/h_{n_k,l})\big)^{1/2}\Big) \\
\le & \mathbb{P}\Big(\max_{n \in N_k} || T_n ||_{\overline{\mathcal{F}}_{k,l,j}} \ge D_1(v)(\tau+\rho)\big(2n_k \log(1/h_{n_k,l})\big)^{1/2}\Big) \\
\le & 4\exp\Big(-A_2 \frac{\rho^2}{\tau^2} \log(1/h_{n_k,l})\Big) \\
= & 4 h_{n_k,l}{}^2.
\end{aligned}$$

Hence, proceeding as in *key argument 2* we get, ultimately as $k \to \infty$,

$$\mathbb{P}_{2,2,k} \le 4 \sum_{0 \le l \le R_k - 1,\ 1 \le j \le J_l} h_{n_k,l}{}^2 \le \frac{4C}{\rho - 1} \mathfrak{h}_{n_{k-1}},$$

whence $(\mathbb{P}_{2,2,k})$ is summable by (HV3).

By making use of similar arguments, it can be proved that, for a suitable choice of $\rho > 1$ and $\delta > 0$ small enough, we have

$$\sum_{k \ge 1} \mathbb{P}_{2,1,k} < \infty. \tag{50}$$

This result is proved by considering the classes

$$\begin{aligned}
\mathcal{F}'_{k,j,l} := & \Big\{ f(z_{k,l,j})^{-1/2} \Big( K\Big(\frac{\cdot - z_{k,l,j}}{h_{n_k,l}{}^{1/d}}\Big) - K\Big(\frac{\cdot - z}{h^{1/d}}\Big)\Big),\ z \in \Gamma_{k,l,j}, \\
& h_{n_k,l} \le h \le \rho h_{n_k,l},\ K \in \mathcal{G} \Big\},
\end{aligned}$$

and showing by $(HV1)$ that, given $\varepsilon > 0$, one can choose $\rho > 1$ and $\delta > 0$ small enough to fulfill $\sup\{\mathrm{Var}\big(g(Z)\big),\ g \in \mathcal{F}'_{k,j,l}\} \leq \epsilon h_{n_k,l}$ uniformly in $j$ and $l$. We omit details for sake of brievness. □

**Remark**:
A close look at the proof of part (i) of Theorem 1 shows that assumption (HK3) can be relaxed to the following assumption:

$$(HK3') \quad \exists M > 0,\ \forall K \in \mathcal{G},\ \forall x \notin [0, M]^d,\ K(x) = 0.$$

### *2.3. Proof of part (ii) of Theorem 1*

Since $\mathbb{K}$ is a compact subset of $\big(\mathcal{B}(\mathcal{G}), ||\cdot||_{\mathcal{G}}\big)$, it is sufficient to show that for fixed $\Psi \in \mathbb{K}$ and $\epsilon > 0$ we have almost surely

$$\lim_{n\to\infty} \sup_{h\in[h_n,\mathfrak{h}_n]} \inf_{z\in H} \left|\left| \frac{G_n(\cdot, h, z)}{\sqrt{2f(z)nh\log(1/h)}} - \Psi \right|\right|_{\mathcal{G}} \leq 4\epsilon.$$

Choose an open hypercube $H' \subset H$ such that $\mathbb{P}(Z_1 \in H') \leq 1/2$. Such a choice is possible because $H$ has a nonempty interior by assumption. Let $1 < \rho$ be a parameter that will be fixed later. Consider the net

$$h_{n,l} := \rho^l h_n,\ l = 0, \ldots, \overline{R}_n - 1,\ \ h_{n,\overline{R}_n} := \mathfrak{h}_n, \tag{51}$$

$$\overline{R}_n := [\log(\mathfrak{h}_n/h_n)/\log(\rho)] + 1. \tag{52}$$

For fixed $l \leq \overline{R}_n$ we divide H' into *disjoint* hypercubes

$$\Gamma_{n,j,l} = z_{n,l,j} + [0, {h_{n,l}}^{1/d})^d.$$

Note than we can construct $J_l := C/h_{n,l}$ disjoint hypercubes, where $C$ depends only on the volume of $H'$.

#### *2.3.1. Step 1*

We shall first show that for any choice of $\rho > 1$ we have almost surely

$$\limsup_{n\to\infty} \sup_{0\leq l\leq \overline{R}_n} \inf_{z\in H} \left|\left| \frac{G_n(\cdot, h_{n,l}, z)}{\big(2f(z)nh_{n,l}\log(1/h_{n,l})\big)^{1/2}} - \Psi \right|\right|_{\mathcal{G}} \leq \epsilon. \tag{53}$$

Recall that $\widetilde{G}_n(\cdot, h, z)$ denotes the "poissonized" version of $G_n(\cdot, h, z)$ (see (38)). By making use of poissonization techniques (see, e.g., Mason (2004), Fact 6), we have, ultimately as $n \to \infty$,

$$\mathbb{P}_n := \mathbb{P}\left( \max_{0\leq l\leq \overline{R}_n} \inf_{z\in H} \left|\left| \frac{G_n(\cdot, h_{n,l}, z)}{\big(2f(z)nh_{n,l}\log(1/h_{n,l})\big)^{1/2}} - \Psi \right|\right|_{\mathcal{G}} > \epsilon \right)$$

$$\leq 2\sum_{l=0}^{\overline{R}_n}\mathbb{P}\Big(\bigcap_{j=1}^{J_l}\Big\|\frac{\widetilde{G}_n(\cdot,h_{n,l},z_{n,l,j})}{\big(2f(z_{n,l,j})nh_{n,l}\log(1/h_{n,l})\big)^{1/2}}-\Psi\Big\|_{\mathcal{G}}>\epsilon\Big)$$
$$=: 2\sum_{l=0}^{\overline{R}_n}\mathbb{P}\Big(\bigcap_{j=1}^{J_l}E_{n,l,j}\Big).$$

But (HK3) entails that, for fixed $l\leq\overline{R}_n,\ j\leq J_l$, the events $E_{n,l,j}$ are mutually independent (by classical properties of Poisson random measures), whence

$$\mathbb{P}\Big(\bigcap_{j=1}^{J_l}E_{n,l,j}\Big)=\prod_{j=1}^{J_l}\Big(1-\mathbb{P}\big(E^C_{n,l,j}\big)\Big)\leq\ \exp\Big(-J_l\min_{j\leq J_l}\mathbb{P}\big(E^C_{n,l,j}\big)\Big).$$

From Proposition 2.3 and by lower semi continuity of $J$ we deduce that, for some $\alpha>0$ and for all large $n$ we have

$$\min_{0\leq l\leq\overline{R}_n,\ 1\leq j\leq J_l}\mathbb{P}\big(E^C_{n,l,j}\big)\geq {h_{n,l}}^{1-\alpha}.$$

Hence, ultimately as $\to\infty$,

$$\begin{aligned}\mathbb{P}_n &\leq \sum_{l=0}^{\overline{R}_n}\exp\big(-J_l {h_{n,l}}^{1-\alpha}\big)\\ &\leq \sum_{l=0}^{\overline{R}_n}\exp\Big(-\Big(\frac{C}{h_{n,l}}-1\Big){h_{n,l}}^{1-\alpha}\Big)\\ &\leq (\overline{R}_n+1)\exp\Big(-\frac{C}{2}{\mathfrak{h}_n}^{-\alpha}\Big)\\ &\leq \Big(1+\frac{\log(\mathfrak{h}_n/h_n)}{\log(\rho)}\Big)\exp\Big(-\frac{C}{2}{\mathfrak{h}_n}^{-\alpha}\Big).\end{aligned}$$

Assumptions (HV2) and (HV3) readily imply that $\mathbb{P}_n$ has a finite sum in $n$, which proves (53) by the Borel-Cantelli lemma.

### *2.3.2. Step 2*

It remains to show that, for any $\epsilon>0$, one can choose $\rho>1$ small enough to have almost surely

$$\overline{\lim}_{n\to\infty}\max_{l\leq\overline{R}_n}\sup_{\substack{z\in H,\\ h_{n,l}\leq h\leq\rho h_{n,l}}}\Big\|\frac{G_n(\cdot,h,z)}{\big(2f(z)nh\log(1/h)\big)^{1/2}}$$
$$-\frac{G_n(\cdot,h_{n,l},z)}{\big(2f(z)nh_{n,l}\log(1/h_{n,l})\big)^{1/2}}\Big\|_{\mathcal{G}}\leq 3\epsilon. \tag{54}$$

We set, for $\rho' \geq 1$ and $x \in \mathbb{R}^d$, $K_{\rho'}(x) := K\big(\rho'^{-1/d}x\big)$. Hence, setting $\rho' = \rho'(h,n) := h/h_{n,l} \in [1,\rho]$ and $u' = u'(h,n) := \big(\log(1/h)\big)^{-1}$, we have almost surely, for each $0 \leq l \leq \overline{R}_n$ and $h \in [h_{n,l}, \rho h_{n,l}]$ and $K \in \mathcal{G}$,

$$\begin{aligned}\frac{G_n(K,h,z)}{\big(2f(z)nh\log(1/h)\big)^{1/2}} &= \left(\frac{h_{n,l}\log(1/h_{n,l})}{h\log(1/h)}\right)^{1/2}\frac{G_n(K_{\rho'},h_{n,l},z)}{\big(2f(z)nh_{n,l}\log(1/h_{n,l})\big)^{1/2}}\\ &= \Big(\rho'^{-1}(1+u'\log(\rho'))\Big)^{1/2}\frac{G_n(K_{\rho'},h_{n,l},z)}{\big(2f(z)nh_{n,l}\log(1/h_{n,l})\big)^{1/2}}\\ &=: \alpha(\rho',u')\frac{G_n(K_{\rho'},h_{n,l},z)}{\big(2f(z)nh_{n,l}\log(1/h_{n,l})\big)^{1/2}}. \qquad (55)\end{aligned}$$

Consider the class

$$\mathcal{G}' := \big\{K_\rho,\ K \in \mathcal{G},\ 1 \leq \rho \leq 2\big\}.$$

Clearly, $\mathcal{G}'$ satisfies (HK1), (HK2), (HK3') and (HK4)–(HK5). By applying part (i) of Theorem 1, we have almost surely

$$\limsup_{n\to\infty}\ \sup_{\substack{h_n\leq h\leq \mathfrak{h}_n,\\ z\in H}}\ \inf_{\Psi\in\mathbb{K}_{\mathcal{G}'}}\Big|\Big|\Psi_n(\cdot,h,z)-\Psi\Big|\Big|_{\mathcal{G}'} = 0, \qquad (56)$$

where $\mathbb{K}_{\mathcal{G}'}$ is as compact subset of $\big(\mathcal{B}(\mathcal{G}'), ||\cdot||_{\mathcal{G}'}\ \big)$ and hence satisfies

$$\lim_{u\to 0,\ \rho\downarrow 1}\ \sup_{\Psi\in\mathbb{K}_{\mathcal{G}'}}\ \sup_{\substack{1\leq\rho'\leq\rho,\\ 0\leq u'\leq u}}\ \sup_{K\in\mathcal{G}} \mid \alpha(\rho',u')\Psi(K_{\rho'}) - \Psi(K) \mid \leq \epsilon. \qquad (57)$$

Now combining (55), (56) and (57) leads to (54), provided that $\rho > 1$ has been chosen small enough. The proof of part (ii) of Theorem 1 is concluded by combining (54) and (53).